\documentclass[a4paper,12pt,reqno]{amsart}
\usepackage[T1]{fontenc}
\usepackage[utf8]{inputenc}
\usepackage[english]{babel}
\usepackage[left=25mm,right=25mm,top=35mm,bottom=35mm]{geometry}
\usepackage{amsmath,amssymb,amsthm}
\usepackage{pgfplots}
\usepackage{pgfplotstable}
\usepackage[margin=0pt,font={small}]{caption}
\usepackage{color,colortbl}
\usepackage[hidelinks]{hyperref}
\usepackage{mathtools}

\usepackage{bigints}
\usepackage{multirow}
\usepackage{mathrsfs}
\usepackage{physics}
\usepackage{tikz}
\usetikzlibrary{arrows}
\usepackage{caption}
\usepackage{subfigure}
\usepackage{stmaryrd}
\usepackage{graphicx} 
\usepackage{booktabs} 

\theoremstyle{definition}
\newtheorem{examp}{Test problem}

\usepackage{fancyhdr}
\advance\footskip0.5cm
\usepackage{color,colortbl}
\definecolor{otherblue}{rgb}{0,0.3,0.6}
\def\rbl#1{\textcolor{otherblue}{#1}}

\def\coloneq{\coloneqq}

\title{IFISS3D:  A computational laboratory for investigating finite element approximation in three dimensions}
\author{Georgios Papanikos}
\address{Department of Mathematics, University of Manchester, Oxford Road, Manchester M13 9PL, UK}
\email{george.papanikos@manchester.ac.uk}

\author{Catherine E. Powell}
\address{Department of Mathematics, University of Manchester, Oxford Road, Manchester M13 9PL, UK}
\email{catherine.powell@manchester.ac.uk}

\author{David J. Silvester}
\address{Department of Mathematics, University of Manchester, Oxford Road, Manchester M13 9PL, UK}
\email{d.silvester@manchester.ac.uk}

\thanks{{\em Acknowledgements.}
This work is supported by   by EPSRC grants
{EP/V048376/1} and {EP/W033801/1}.
}

\date{\today}

\begin{document}

\begin{abstract}
IFISS is an established MATLAB finite element software package for studying strategies for solving partial differential equations (PDEs). IFISS3D is a new add-on toolbox that extends IFISS capabilities for elliptic PDEs from two to three space dimensions. The open-source MATLAB framework provides a computational laboratory for experimentation and exploration of finite element approximation and error estimation, as well as iterative solvers.  The package is designed to be useful as a teaching tool for instructors and students who want to learn about state-of-the-art finite element methodology.  It will also be useful for researchers as a source of reproducible test matrices of arbitrarily large dimension.
\end{abstract}

\maketitle
\thispagestyle{fancy}

\section{Introduction and brief history}

The IFISS software \cite{ifiss} was developed by Elman, Ramage and Silvester~\cite{ers07}.
 It can  be run in MATLAB (developed by the MathWorks${}^\copyright$) or  Gnu Octave (free software).
It is  structured as a stand-alone package for studying discretisation 
algorithms for partial differential equations (PDEs), and for exploring and developing algorithms in
numerical linear and nonlinear algebra for solving the associated discrete systems. It can be used
as a pedagogical tool for studying these issues, or more elementary ones such as the properties of Krylov 
subspace iterative methods.   Investigative numerical experiments in a teaching setting
enable students to develop deduction and interpretation skills, and are especially useful in helping 
students  to remember critical ideas in the  long term.
 IFISS is also an established starting point for developing code for specialised research applications
 (as evidenced by the variety of citations to it, see~\cite{swmath}), and is extensively used 
 by researchers  in numerical linear algebra as a source of reproducible test matrices of arbitrarily large 
 dimension.

The development of the MATLAB functionality during the period 1990--2005 opened up the possibility of creating a problem-based-learning  environment (notably the  IFISS package)
 that could be used together with standard teaching mechanisms to facilitate understanding of abstract 
 theoretical concepts. 
The functionality of IFISS  was significantly  extended  in the  
 period between 2005 and 2015---culminating in the publication of the review article~\cite{ers14},  
 which coincided with the publication of  the {second edition} of the monograph~\cite{elman15}. 

A unique feature of IFISS is its comprehensive  nature. For each problem it addresses, 
 it enables the study of both discretisation and iterative solution algorithms, as well as the interaction
between the two and the resulting effect on solution cost. However, it is restricted
 to the solution of PDEs on two-dimensional spatial domains.  This limitation can be overcome by adding the new IFISS3D toolbox \cite{ifiss3D} to the existing IFISS software. The three-dimensional finite element approximation and error estimation strategies included in the new software are specified in the next section. Section~\ref{sec.referenceproblems}
 describes three reference problems that provide a convenient
 starting point for studying  rates of convergence of the approximations  to the true solution. 
 The  structure of the IFISS3D package is discussed  in Section~\ref{sec.structure}. The directory
 structure is intended to simplify the task of extending the functionality to other PDE problems 
 and higher-order finite element methods.  Case studies of two important aspects of 
 three-dimensional finite element approximation are presented in Section~\ref{sec.casestudies}. 
 
\section{Discretisation and Error Estimation Specifics} \label{sec.specifics}
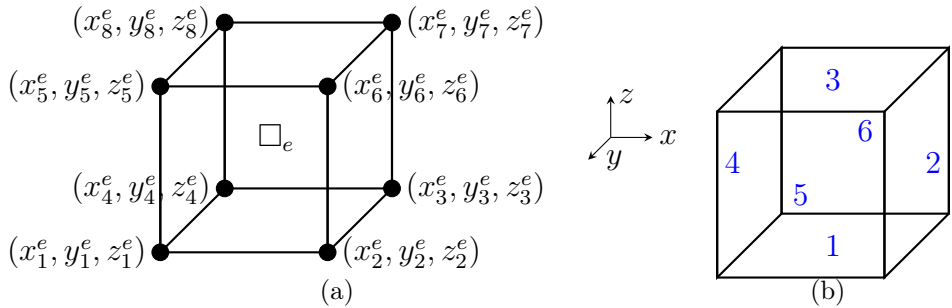
\begin{figure}[!ht]
\begin{tabular}{  c c }
		\subfigure[]{\begin{tikzpicture}[scale=1.1]
		\draw[thick](-1,-1,-1)--(1,-1,-1)--(1,1,-1)--(-1,1,-1)--(-1,-1,-1);
		\draw[thick] (-1,-1,1)--(1,-1,1)--(1,1,1)--(-1,1,1)--(-1,-1,1);
		\draw[thick] (-1,-1,-1)--(-1,1,-1)--(-1,1,1)--(-1,-1,1)--(-1,-1,-1);
		\draw[thick] (1,-1,-1)--(1,1,-1)--(1,1,1)--(1,-1,1)--(1,-1,-1);		
		\draw(0,0,0) node{$\square_{e}$};
		\draw(-2,-1,1) node{ $(x^e_1,y^e_1,z^e_1)$ };
		\draw(2,-1,1) node{$(x^e_2,y^e_2,z^e_2)$};
		\draw(2,-1,-1) node{$(x^e_3,y^e_3,z^e_3)$};
		\draw(-2,-1,-1) node{$(x^e_4,y^e_4,z^e_4)$};
		\draw(-2,1,1) node{$(x^e_5,y^e_5,z^e_5)$};
		\draw(2,1,1) node{$(x^e_6,y^e_6,z^e_6)$};
		\draw(-2,1,-1) node{$(x^e_8,y^e_8,z^e_8)$};
		\draw(2,1,-1) node{$(x^e_7,y^e_7,z^e_7)$};	
		\filldraw[color=black](-1,-1,1)  circle (0.1cm);
		\filldraw[color=black](1,-1,1) circle (0.1cm);
		\filldraw[color=black](1,-1,-1) circle (0.1cm);
	 	\filldraw[color=black](-1,-1,-1) circle (0.1cm);
		\filldraw[color=black](-1,1,1)  circle (0.1cm);
          	\filldraw[color=black](1,1,1) circle (0.1cm);
		\filldraw[color=black](-1,1,-1) circle (0.1cm);
		\filldraw[color=black](1,1,-1) circle (0.1cm);
		\draw [->,>=stealth] (4,0,0) -- (4.5,0,0) node{$ \ \  \ x$};
		\draw [->,>=stealth] (4,0,0) -- (4,0.5,0) node{$ \ \ \ z$};
		\draw [->,>=stealth] (4,0,0) -- (4,0,0.7) node{$ \ \ \ \  \ y$};	
		\end{tikzpicture}}	
			& 	 \subfigure[]{\begin{tikzpicture}[scale=1.1]
			\draw[thick](-1,-1,-1)--(1,-1,-1)--(1,1,-1)--(-1,1,-1)--(-1,-1,-1);
			\draw[thick] (-1,-1,1)--(1,-1,1)--(1,1,1)--(-1,1,1)--(-1,-1,1);
			\draw[thick] (-1,-1,-1)--(-1,1,-1)--(-1,1,1)--(-1,-1,1)--(-1,-1,-1);
			\draw[thick] (1,-1,-1)--(1,1,-1)--(1,1,1)--(1,-1,1)--(1,-1,-1);
			\draw[blue](0,0,1) node{5};
			\draw[blue](0,1,0) node{3};
			\draw[blue](1.2,0,0) node{2};
			\draw[blue](0,-1,0) node{1};
			\draw[blue](-1.2,0,0) node{4};
			\draw[blue](0,0,-1) node{6};
			\end{tikzpicture}}  	
		\end{tabular}
		\vspace{-2mm}
	\caption{ Hexahedral (brick) element (a) vertex numbering and (b) face numbering.}
	\label{fig:elementnumbering}
\end{figure}

The IFISS3D software generates  approximations to the solution of PDEs
modelling physical problems in three spatial dimensions. The starting point for the  process
is a finite element partitioning  of a domain of interest  $D\subset\mathbb{R}^3$ into $n_e$ 
hexahedral (brick) elements $\square_e\subset D$, $e=1,2,3, \ldots,n_e$ so that
	\begin{equation}\nonumber
		\begin{cases}
			\square_e \ \text{is open in} \ \mathbb{R}^3 \\
			D = \overline{\cup_{e=1}^{n_e}\square_e} \\ 
			\square_{i}\cap \square_{j} = \emptyset, \ \ i\neq j,
		\end{cases}
	\end{equation}
where the upper bar represents the closure of the union. 
An arbitrary  element $\square_{e}$ is a hexahedron with six faces and with local vertex coordinates 
$(x^e_i,y^e_i,z^e_i)$, $i=1,2, \ldots, 8$ ordered as shown in Fig.~\ref{fig:elementnumbering}.

The simplest choice of a conforming finite element space in $\mathbb{R}^3$ is the 
$\mathbb{Q}_1$ approximation space of  piecewise  {\it trilinear} polynomials that  
take the form
	\begin{equation} \nonumber 
	q_1^e (\boldsymbol{x}) =	(\alpha_1x+\alpha_2)(\alpha_3y+\alpha_4)(\alpha_5z+\alpha_6)
	\end{equation}	
on each element $\square_{e}$. The continuity of the global approximation  is ensured by  
defining a Lagrangian basis for $q^e_1$  at the eight vertices  of the hexahedron, that is
$$ \phi^e_i(\boldsymbol{x}) = \begin{cases}
		1 & \text{if} \ \ x = x^e_i, y = y^e_i, \ z= z^e_i, \\
		0 & \text{at the other vertices}
	\end{cases}, \quad i =1,2,\ldots, 8 .
$$

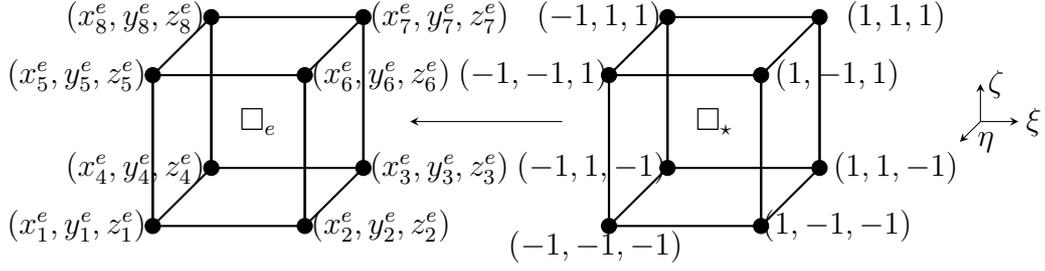
\begin{figure}	
	\begin{center}
		\begin{tikzpicture}
		\draw[thick](-1,-1,-1)--(1,-1,-1)--(1,1,-1)--(-1,1,-1)--(-1,-1,-1);
		\draw[thick] (-1,-1,1)--(1,-1,1)--(1,1,1)--(-1,1,1)--(-1,-1,1);
		\draw[thick] (-1,-1,-1)--(-1,1,-1)--(-1,1,1)--(-1,-1,1)--(-1,-1,-1);
		\draw[thick] (1,-1,-1)--(1,1,-1)--(1,1,1)--(1,-1,1)--(1,-1,-1);		
		\draw(0,0,0) node{$\square_{e}$};
		\draw(-2,-1,1) node{ $(x^e_1,y^e_1,z^e_1)$ };
		\draw(2,-1,1) node{$(x^e_2,y^e_2,z^e_2)$};
		\draw(2,-1,-1) node{$(x^e_3,y^e_3,z^e_3)$};
		\draw(-2,-1,-1) node{$(x^e_4,y^e_4,z^e_4)$};
		\draw(-2,1,1) node{$(x^e_5,y^e_5,z^e_5)$};
		\draw(2,1,1) node{$(x^e_6,y^e_6,z^e_6)$};
		\draw(-2,1,-1) node{$(x^e_8,y^e_8,z^e_8)$};
		\draw(2,1,-1) node{$(x^e_7,y^e_7,z^e_7)$};
			
		\filldraw[color=black](-1,-1,1)  circle (0.1cm);		
		\filldraw[color=black](1,-1,1) circle (0.1cm);
		\filldraw[color=black](1,-1,-1) circle (0.1cm);	
		\filldraw[color=black](-1,-1,-1) circle (0.1cm);
		\filldraw[color=black](-1,1,1)  circle (0.1cm);	
		\filldraw[color=black](1,1,1) circle (0.1cm);	
		\filldraw[color=black](-1,1,-1) circle (0.1cm);		
		\filldraw[color=black](1,1,-1) circle (0.1cm);
				
			
			\draw[thick] (5,-1,-1)--(7,-1,-1)--(7,1,-1)--(5,1,-1)--(5,-1,-1);
			\draw[thick] (5,-1,1)--(7,-1,1)--(7,1,1)--(5,1,1)--(5,-1,1);
			\draw[thick] (5,-1,-1)--(5,1,-1)--(5,1,1)--(5,-1,1)--(5,-1,-1);
			\draw[thick] (7,-1,-1)--(7,1,-1)--(7,1,1)--(7,-1,1)--(7,-1,-1);
						
			\filldraw[color=black](5,-1,1)  circle (0.1cm);	
			\filldraw[color=black](7,-1,1) circle (0.1cm);	
			\filldraw[color=black](7,-1,-1) circle (0.1cm);	
			\filldraw[color=black](5,-1,-1) circle (0.1cm);	
			\filldraw[color=black](5,1,1)  circle (0.1cm);	
			\filldraw[color=black](7,1,1) circle (0.1cm);
			\filldraw[color=black](5,1,-1) circle (0.1cm);
			\filldraw[color=black](7,1,-1) circle (0.1cm);
			
			\draw(6,0,0) node{$\square_{\star}$};
			 (0.1cm);\filldraw[color=black](4.9,-1.2,1.2)  circle node{ $(-1,-1,-1)$ };
			\draw(8,-1,1) node{$(1,-1,-1)$};
			\draw(8,-1,-1) node{$(1,1,-1)$};
			\draw(4,-1,-1) node{$(-1,1,-1)$};
			\draw(4,1,1) node{$(-1,-1,1)$};
			\draw(8,1,1) node{$(1,-1,1)$};
			\draw(4.1,1,-1) node{$(-1,1,1)$};
			\draw(8,1,-1) node{$(1,1,1)$};				
			\draw [<-,>=stealth] (2,0,0) -- (4,0,0);
			\draw [->,>=stealth] (9.5,0,0) -- (10,0,0) node{$ \ \ \ \xi$};
			\draw [->,>=stealth] (9.5,0,0) -- (9.5,0.5,0) node{$ \ \ \ \zeta$};
			\draw [->,>=stealth] (9.5,0,0) -- (9.5,0,0.7) node{$ \ \ \ \ \ \eta$};			
		\end{tikzpicture}
	\end{center} 	
	\vspace{-2mm}
	\caption{Isoparametric mapping from a reference cube.}\label{fig.isomapping}	
 \end{figure}

	Additional geometric flexibility ({\it stretched} grids) can be incorporated by constructing an {\it isoparametric} transformation 
from the reference cube  $[-1,1]^3$ (denoted  $\square_{\star}$)  to a general hexahedron  $\square_{e}$.  To this end we define the following basis functions for the reference element 
\begin{align}
\ell^i(\xi,\eta,\zeta) = 1/8 (1+\xi_i\xi)&(1+\eta_i\eta)(1+\zeta_i\zeta), \quad i=1,2,\ldots ,8 \nonumber
\end{align} 
where $\xi_i,  \eta_i,  \zeta_i$ are node values of $\square_{\star}$ and $\xi,\eta,\zeta\in [-1,1]$ and map to an arbitrary hexahedral element $\square_{e}$  with 
vertices  $(x^e_i,y^e_i,z^e_i)$, $i=1,2,\ldots, 8$ by the change of variables	
	\[
	x(\xi,\eta,\zeta) = \sum_{i=1}^8 x^e_i \ell^i(\xi,\eta,\zeta),  \ y(\xi,\eta,\zeta) 
	= \sum_{i=1}^8 y^e_i \ell^i(\xi,\eta,\zeta), \  z(\xi,\eta,\zeta) 
	= \sum_{i=1}^8 z^e_i \ell^i(\xi,\eta,\zeta).
	\] 
This mapping is  illustrated in Fig.~\ref{fig.isomapping}.

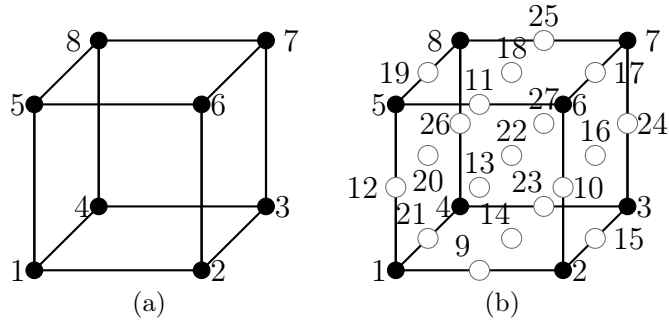
\begin{figure}[!th]
		\centering
		\begin{tabular}{  c c   }
			 \subfigure[]{\begin{tikzpicture}[scale=1.1]
			
				\draw[thick](-1,-1,-1)--(1,-1,-1)--(1,1,-1)--(-1,1,-1)--(-1,-1,-1);
				\draw[thick] (-1,-1,1)--(1,-1,1)--(1,1,1)--(-1,1,1)--(-1,-1,1);
				\draw[thick] (-1,-1,-1)--(-1,1,-1)--(-1,1,1)--(-1,-1,1)--(-1,-1,-1);
				\draw[thick] (1,-1,-1)--(1,1,-1)--(1,1,1)--(1,-1,1)--(1,-1,-1);
				
				\filldraw[color=black](-1,-1,1)  circle (0.1cm);
				\filldraw[color=black](-1.2,-1,1) node{1};	
				\filldraw[color=black](1,-1,1) circle (0.1cm);
				\filldraw[color=black](1.2,-1,1) node{2};	
				\filldraw[color=black](1,-1,-1) circle (0.1cm);
				\filldraw[color=black](1.2,-1,-1) node{3};	
				\filldraw[color=black](-1,-1,-1) circle (0.1cm);
				\filldraw[color=black](-1.2,-1,-1) node{4};	
				\filldraw[color=black](-1,1,1)  circle (0.1cm);
				\filldraw[color=black](-1.2,1,1) node{5};	
				\filldraw[color=black](1,1,1) circle (0.1cm);
				\filldraw[color=black](1.2,1,1) node{6};	
				\filldraw[color=black](-1,1,-1) circle (0.1cm);
				\filldraw[color=black](-1.3,1,-1) node{8};	
				\filldraw[color=black](1,1,-1) circle (0.1cm);
				\filldraw[color=black](1.3,1,-1) node{7};
				
			\end{tikzpicture}}	
			& 	 \subfigure[]{\begin{tikzpicture}[scale=1.1]
				
				\draw[thick](-1,-1,-1)--(1,-1,-1)--(1,1,-1)--(-1,1,-1)--(-1,-1,-1);
				\draw[thick] (-1,-1,1)--(1,-1,1)--(1,1,1)--(-1,1,1)--(-1,-1,1);
				\draw[thick] (-1,-1,-1)--(-1,1,-1)--(-1,1,1)--(-1,-1,1)--(-1,-1,-1);
				\draw[thick] (1,-1,-1)--(1,1,-1)--(1,1,1)--(1,-1,1)--(1,-1,-1);
				
				\filldraw[color=black](-1,-1,1)  circle (0.1cm);
				\filldraw[color=black](-1.2,-1,1) node{1};	
				\filldraw[color=black](1,-1,1) circle (0.1cm);
				\filldraw[color=black](1.2,-1,1) node{2};	
				\filldraw[color=black](1,-1,-1) circle (0.1cm);
				\filldraw[color=black](1.2,-1,-1) node{3};	
				\filldraw[color=black](-1,-1,-1) circle (0.1cm);
				\filldraw[color=black](-1.2,-1,-1) node{4};	
				\filldraw[color=black](-1,1,1)  circle (0.1cm);
				\filldraw[color=black](-1.2,1,1) node{5};	
				\filldraw[color=black](1,1,1) circle (0.1cm);
				\filldraw[color=black](1.2,1,1) node{6};		
				\filldraw[color=black](-1,1,-1) circle (0.1cm);
				\filldraw[color=black](-1.3,1,-1) node{8};	
				\filldraw[color=black](1,1,-1) circle (0.1cm);
				\filldraw[color=black](1.3,1,-1) node{7};	
				\filldraw[color=black!60, fill=black!-100, very thin](0,-1,1) circle (0.12cm);
				\filldraw[color=black](-0.2,-0.7,1) node{9};			
				\filldraw[color=black!60, fill=black!-100, very thin](1,0,1)  circle (0.12cm);
				\filldraw[color=black](1.3,0,1) node{10};	
				\filldraw[color=black!60, fill=black!-100, very thin](0,1,1)  circle (0.12cm);
				\filldraw[color=black](0,1.3,1) node{11};	
				\filldraw[color=black!60, fill=black!-100, very thin](-1,0,1)  circle (0.12cm);
				\filldraw[color=black](-1.4,0,1) node{12};	
				\filldraw[color=black!60, fill=black!-100, very thin](0,0,1) circle (0.12cm);
				\filldraw[color=black](0,0.3,1) node{13};	
				\filldraw[color=black!60, fill=black!-100, very thin](0,-1,0) circle (0.12cm);
				\filldraw[color=black](-0.2,-0.7,0) node{14};	
				\filldraw[color=black!60, fill=black!-100, very thin](1,-1,0)  circle (0.12cm);
				\filldraw[color=black](1.4,-1,0) node{15};	
				\filldraw[color=black!60, fill=black!-100, very thin](1,0,0)  circle (0.12cm);
				\filldraw[color=black](1,0.3,0) node{16};	
				\filldraw[color=black!60, fill=black!-100, very thin](1,1,0)  circle (0.12cm);
				\filldraw[color=black](1.4,1,0) node{17};	
				\filldraw[color=black!60, fill=black!-100, very thin](0,1,0)  circle (0.12cm);
				\filldraw[color=black](0,1.3,0) node{18};	
				\filldraw[color=black!60, fill=black!-100, very thin](-1,1,0)  circle (0.12cm);
				\filldraw[color=black](-1.4,1,0) node{19};	
				\filldraw[color=black!60, fill=black!-100, very thin](-1,0,0)  circle (0.12cm);
				\filldraw[color=black](-1,-0.3,0) node{20};	
				\filldraw[color=black!60, fill=black!-100, very thin](-1,-1,0)  circle (0.12cm);
				\filldraw[color=black](-1.2,-0.7,0) node{21};	
				\filldraw[color=black!60, fill=black!-100, very thin](0,0,0) circle (0.12cm);
				\filldraw[color=black](0,0.3,0) node{22};	
				\filldraw[color=black!60, fill=black!-100, very thin](0,-1,-1) circle (0.12cm);
				\filldraw[color=black](-0.2,-0.7,-1) node{23};	
				\filldraw[color=black!60, fill=black!-100, very thin](1,0,-1)  circle (0.12cm);
				\filldraw[color=black](1.3,0,-1) node{24};	
				\filldraw[color=black!60, fill=black!-100, very thin](0,1,-1)  circle (0.12cm);
				\filldraw[color=black](0,1.3,-1) node{25};
				\filldraw[color=black!60, fill=black!-100, very thin](-1,0,-1)  circle (0.12cm);
				\filldraw[color=black](-1.3,0,-1) node{26};
				\filldraw[color=black!60, fill=black!-100, very thin](0,0,-1) circle (0.12cm);
				\filldraw[color=black](0,0.3,-1) node{27};
			\end{tikzpicture}}  		
		\end{tabular}
		\vspace{-2mm}
	\caption{ (a)  $\mathbb{Q}_1$ nodes and numbering and (b) $\mathbb{Q}_2$ nodes and numbering.}
	\label{fig:Q1Q2nodes}
\end{figure}

The IFISS3D software also provides the option of  higher-order approximation using
polynomials that are globally continuous but piecewise {\it triquadratic}  ($\mathbb{Q}_2$), that is,
\begin{equation}\nonumber 
 q_2^e (\boldsymbol{x}) =  (\alpha_1x^2 + \alpha_2x + \alpha_3)
                                          (\alpha_4y^2 + \alpha_5y + \alpha_6)(\alpha_7z^2 + \alpha_8z + \alpha_9).
\end{equation}
The continuity of the global  $\mathbb{Q}_2$ approximation  is ensured by  
defining a Lagrangian basis for $q_2^e$  at the eight vertices  of the hexahedron together  with
the nineteen additional nodes shown in Fig.~\ref{fig:Q1Q2nodes}.

The  $\mathbb{Q}_2$  isoparametric transformation is given  by
$$
	x(\xi,\eta,\zeta) = \sum_{i=1}^{27} x^e_i {\psi}^i(\xi,\eta,\zeta),  
	\quad y(\xi,\eta,\zeta) = \sum_{i=1}^{27 }y^e_i {\psi}^i(\xi,\eta,\zeta), 
	\quad z(\xi,\eta,\zeta) = \sum_{i=1}^{27} z^e_i {\psi}^i(\xi,\eta,\zeta), 
$$ 
with reference basis functions
\begin{align*}
{\psi}^i(\xi,\eta,\zeta) = {N}^k(\xi) \, {N}^l(\eta) \, {N}^m(\zeta),\ i=1,2,\ldots,27; \  k,l,m =1,2,3  
\end{align*}
where
$$
{N}^1(\xi) = 1/2(\xi-1)\xi, \ {N}^2(\xi) = 1-\xi^2, \ \ {N}^3(\xi) = 1/2(1+\xi)\xi.
$$
	
A fundamental feature of  the IFISS software is the use of hierarchical error estimation. This strategy  was 
developed for scalar elliptic PDEs by Bank \& Smith~\cite{bank93}, but has been extended
to more general PDE problems (including systems of PDEs) over the past two decades. Crucially, the hierarchical approach yields reliable estimates of the {\it error reduction} that
can be expected using an enhanced approximation. It also provides a rigorous setting for establishing
the convergence of adaptive refinement strategies, such as those that are built into the T-IFISS package~\cite{bespalov20}, 
and the ML-SGFEM software~\cite{multilevelsg} associated 
 with the work of Crowder et al.~\cite{crowder19} on multilevel stochastic Galerkin finite element methods for parametric PDEs.
 
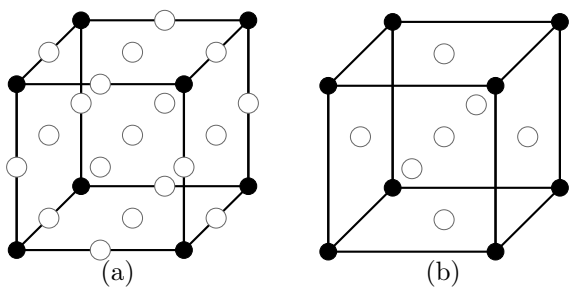
\begin{figure}[!ht]
	\centering
	\begin{tabular}{c c}	
   	\subfigure[]{
   	\begin{tikzpicture}[scale=1.1]
   		\draw[thick](-1,-1,-1)--(1,-1,-1)--(1,1,-1)--(-1,1,-1)--(-1,-1,-1);
   		\draw[thick] (-1,-1,1)--(1,-1,1)--(1,1,1)--(-1,1,1)--(-1,-1,1);
   		\draw[thick] (-1,-1,-1)--(-1,1,-1)--(-1,1,1)--(-1,-1,1)--(-1,-1,-1);
   		\draw[thick] (1,-1,-1)--(1,1,-1)--(1,1,1)--(1,-1,1)--(1,-1,-1);
   		
   		\filldraw[color=black](-1,-1,1)  circle (0.1cm);
   		\filldraw[color=black](-1.2,-1,1);	
   		\filldraw[color=black](1,-1,1) circle (0.1cm);
   		\filldraw[color=black](1.2,-1,1);	
   		\filldraw[color=black](1,-1,-1) circle (0.1cm);
   		\filldraw[color=black](1.2,-1,-1);	
   		\filldraw[color=black](-1,-1,-1) circle (0.1cm);
   		\filldraw[color=black](-1.2,-1,-1);
   		\filldraw[color=black](-1,1,1)  circle (0.1cm);
   		\filldraw[color=black](-1.2,1,1);
   		\filldraw[color=black](1,1,1) circle (0.1cm);
   		\filldraw[color=black](1.2,1,1);
   		\filldraw[color=black](-1,1,-1) circle (0.1cm);
   		\filldraw[color=black](-1.3,1,-1);
   		\filldraw[color=black](1,1,-1) circle (0.1cm);
   		\filldraw[color=black](1.3,1,-1) ;
   		
   		\filldraw[color=black!60, fill=black!-100, very thin](0,-1,1) circle (0.12cm);
   		\filldraw[color=black](-0.2,-0.7,1) ;
   		\filldraw[color=black!60, fill=black!-100, very thin](1,0,1)  circle (0.12cm);
   		\filldraw[color=black](1.3,0,1) ;
   		\filldraw[color=black!60, fill=black!-100, very thin](0,1,1)  circle (0.12cm);
   		\filldraw[color=black](0,1.3,1);
   		\filldraw[color=black!60, fill=black!-100, very thin](-1,0,1)  circle (0.12cm);
   		\filldraw[color=black](-1.4,0,1);
   		\filldraw[color=black!60, fill=black!-100, very thin](0,0,1) circle (0.12cm);
   		\filldraw[color=black](0,0.3,1);
   		\filldraw[color=black!60, fill=black!-100, very thin](0,-1,0) circle (0.12cm);
   		\filldraw[color=black](-0.2,-0.7,0);
   		\filldraw[color=black!60, fill=black!-100, very thin](1,-1,0)  circle (0.12cm);
   		\filldraw[color=black](1.4,-1,0) ;
   		\filldraw[color=black!60, fill=black!-100, very thin](1,0,0)  circle (0.12cm);
   		\filldraw[color=black](1,0.3,0);
   		\filldraw[color=black!60, fill=black!-100, very thin](1,1,0)  circle (0.12cm);
   		\filldraw[color=black](1.4,1,0) ;
   		\filldraw[color=black!60, fill=black!-100, very thin](0,1,0)  circle (0.12cm);
   		\filldraw[color=black](0,1.3,0);
   		\filldraw[color=black!60, fill=black!-100, very thin](-1,1,0)  circle (0.12cm);
   		\filldraw[color=black](-1.4,1,0) ;
   		\filldraw[color=black!60, fill=black!-100, very thin](-1,0,0)  circle (0.12cm);
   		\filldraw[color=black](-1,-0.3,0) ;
   		\filldraw[color=black!60, fill=black!-100, very thin](-1,-1,0)  circle (0.12cm);
   		\filldraw[color=black](-1.2,-0.7,0) ;
   		\filldraw[color=black!60, fill=black!-100, very thin](0,0,0) circle (0.12cm);
   		\filldraw[color=black](0,0.3,0);
   		\filldraw[color=black!60, fill=black!-100, very thin](0,-1,-1) circle (0.12cm);
   		\filldraw[color=black](-0.2,-0.7,-1);
   		\filldraw[color=black!60, fill=black!-100, very thin](1,0,-1)  circle (0.12cm);
   		\filldraw[color=black](1.3,0,-1) ;
   		\filldraw[color=black!60, fill=black!-100, very thin](0,1,-1)  circle (0.12cm);
   		\filldraw[color=black](0,1.3,-1);
   		\filldraw[color=black!60, fill=black!-100, very thin](-1,0,-1)  circle (0.12cm);
   		\filldraw[color=black](-1.3,0,-1);
   		\filldraw[color=black!60, fill=black!-100, very thin](0,0,-1) circle (0.12cm);
   		\filldraw[color=black](0,0.3,-1) ;
   	\end{tikzpicture}} 
   	
 & \subfigure[]{
   	\begin{tikzpicture}[scale=1.1]
   		
   		\draw[thick](-1,-1,-1)--(1,-1,-1)--(1,1,-1)--(-1,1,-1)--(-1,-1,-1);
   		\draw[thick] (-1,-1,1)--(1,-1,1)--(1,1,1)--(-1,1,1)--(-1,-1,1);
   		\draw[thick] (-1,-1,-1)--(-1,1,-1)--(-1,1,1)--(-1,-1,1)--(-1,-1,-1);
   		\draw[thick] (1,-1,-1)--(1,1,-1)--(1,1,1)--(1,-1,1)--(1,-1,-1);
   		
   		\filldraw[color=black](-1,-1,1)  circle (0.1cm);
   		\filldraw[color=black](-1.2,-1,1) ;
   		\filldraw[color=black](1,-1,1) circle (0.1cm);
   		\filldraw[color=black](1.2,-1,1);
   		\filldraw[color=black](1,-1,-1) circle (0.1cm);
   		\filldraw[color=black](1.2,-1,-1);
   		\filldraw[color=black](-1,-1,-1) circle (0.1cm);
   		\filldraw[color=black](-1.2,-1,-1);
   		\filldraw[color=black](-1,1,1)  circle (0.1cm);
   		\filldraw[color=black](-1.2,1,1) ;
   		\filldraw[color=black](1,1,1) circle (0.1cm);
   		\filldraw[color=black](1.2,1,1) ;
   		\filldraw[color=black](-1,1,-1) circle (0.1cm);
   		\filldraw[color=black](-1.3,1,-1) ;
   		\filldraw[color=black](1,1,-1) circle (0.1cm);
   		\filldraw[color=black](1.3,1,-1) ;
   		
   		\filldraw[color=black!60, fill=black!-100, very thin](0,0,1) circle (0.12cm);
   		\filldraw[color=black](0,0.3,1);
   		\filldraw[color=black!60, fill=black!-100, very thin](0,-1,0) circle (0.12cm);
   		\filldraw[color=black](-0.2,-0.7,0);
   		\filldraw[color=black!60, fill=black!-100, very thin](1,0,0)  circle (0.12cm);
   		\filldraw[color=black](1,0.3,0);
   		\filldraw[color=black!60, fill=black!-100, very thin](0,1,0)  circle (0.12cm);
   		\filldraw[color=black](0,1.3,0);
   		\filldraw[color=black!60, fill=black!-100, very thin](-1,0,0)  circle (0.12cm);
   		\filldraw[color=black](-1,-0.3,0);
   		\filldraw[color=black!60, fill=black!-100, very thin](0,0,0) circle (0.12cm);
   		\filldraw[color=black](0,0.3,0);
   		\filldraw[color=black!60, fill=black!-100, very thin](0,0,-1) circle (0.12cm);
   		\filldraw[color=black](0,0.3,-1);
   	\end{tikzpicture} }
\end{tabular}
\vspace{-2mm}
\caption{White nodes associated with bubble functions used to construct (a) full $\mathbb{Q}_2(h)$ or $\mathbb{Q}_1(h/2)$ and (b) reduced $\mathbb{Q}^r_2(h)$ and $\mathbb{Q}^r_1(h/2)$ error estimation spaces.}	
\label{fig:spacesYk}
\end{figure}

A posteriori error estimation in  practical finite element software 
(such as DUNE~\cite{dune16} or FEniCS~\cite{logg12}) is typically done using  residual 
error estimation strategies. This requires the computation of  norms of PDE residuals in the interior of each element and norms of flux jumps 
(edge residuals) on inter-element faces.  The additional computational cost of  hierachical error estimation is nontrivial.
Having generated a solution using $\mathbb{Q}_1$ approximation\footnote{Hierachical error estimation for
 $\mathbb{Q}_2$ approximation is not included in the current release of IFISS3D.},  computed
interior and edge residuals  are  input as source data for  {\it element} PDE problems  
that are solved numerically using  an enhanced approximation space.  In IFISS3D, one can construct the enhanced space using triquadratic  basis functions on the original element
 $(\mathbb{Q}_2(h))$,  or trilinear basis functions defined on a subdivision of the original element into 8 smaller ones $(\mathbb{Q}_1(h/2))$.
 These basis or `bubble' functions are associated with the white nodes illustrated in Fig.~\ref{fig:spacesYk}(a), leading to linear algebra systems of dimension nineteen. Alternatively, {\it reduced} versions of these spaces of dimension 7, denoted $\mathbb{Q}^r_2(h)$ and $\mathbb{Q}^r_1(h/2)$, can also be constructed by incorporating only the basis functions associated with the interior node and 
the central nodes on each face, as illustrated in Fig.~\ref{fig:spacesYk}(b). In all four cases, a low-dimensional system
must be solved for every element in the mesh. This calculation is efficiently vectorised in IFISS. 

The computational overhead of hierarchical error estimation will be discussed further in Section~4.

\section{Reference problems}\label{sec.referenceproblems}

Three test problems are built into the IFISS3D toolbox.  Illustrative results for these problems are discussed in this section. 
 The reported timings were obtained on a 2.9 GHz 6-Core Intel Core i9  
MacBook  using the {\tt tic} {\tt toc} functionality built into MATLAB. 
In all cases, we compute an approximation $u_h \in X_h$ to $u\in X \subset H^1(D)$ satisfying
the  standard weak formulation of  the following Poisson problem
	\begin{align}\label{poisson}
		-\nabla^2 u &= 1  \quad\text{in } D\subset\mathbb{R}^3
			\\
			u &= 0  \quad \text{on } \partial D. \label{bc}
	\end{align}

\begin{examp}[convex domain]\label{ex.cubedomain}
A  $\mathbb{Q}_1$ finite element solution to (\ref{poisson})--(\ref{bc})  defined  on 
the domain  $D=[-1,1]^3$  
is shown in Fig.~\ref{fig.testproblem1}. In this computation
 the cube-shaped domain has been subdivided uniformly into  $32^3$ elements. The   dimension
 of the resulting linear algebra  system  is  35,937 (the boundary vertices are retained when
 assembling the system).  The MATLAB R2021b sparse direct solver ($\backslash$)  solves this  
 system in about half a second. The solution and estimated errors plotted in 
 the cross section shown in  Fig.~\ref{fig.testproblem1}  are consistent with the plots that are generated 
  when solving   (\ref{poisson})--(\ref{bc})   using IFISS software on the two-dimensional domain $D=[-1,1]\times[-1,1]$; 
 see Fig~1.1 in \cite{elman15}. 
 As expected, the largest error is concentrated near the sides of the cube.

 \begin{figure}[!ht]
       \begin{center}
	\includegraphics[width=0.85\linewidth]{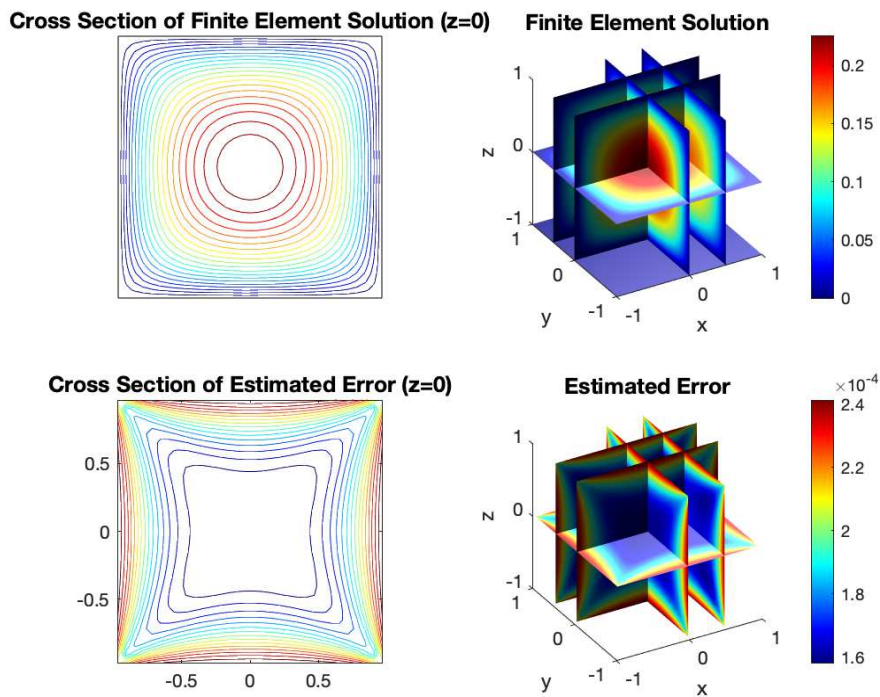}		 
	\end{center}\vspace{-4mm}
\caption{Solution and estimated energy error distribution for test problem~1.}
\label{fig.testproblem1}
\end{figure}
 
  Exploiting  Galerkin orthogonality, the exact  energy 
error can be estimated  by comparing the energy of a reference solution\footnote{ 
$ \|u_{\rm ref}\|^2_E = 0.64539192 $   computed using $\mathbb{Q}_2$  
approximation on  a grid of 64$^3$ elements.}
with the energy of the computed finite element solution
$
 \| u_{\rm ref}-u_h \|_E^2 =  \| u_{\rm ref}\|_E^2 -  \| u_h \|_E^2.
$
Estimated errors  that are computed using this strategy
 are presented in Table~\ref{tab.testproblem1errors}.

 \begin{table}[!ht]
 \caption{Estimated energy errors for convex domain test problem.}
 \vspace{-2mm}
 \label{tab.testproblem1errors}
	\begin{center}
		\begin{tabular}{  c c c c | c c c} 
		& & \multicolumn{2}{c}{$\mathbb{Q}_1\qquad$} 
		&  \multicolumn{3}{c}{$\mathbb{Q}_2 $}  \\
		$n_e$ & $h$ & 
		                    $\|u_h\|^2_E$ & $\| u_{\rm ref}-u_h \|_E$  &
		      $n_e$ & $\|u_h\|^2_E$ & $\| u_{\rm ref}-u_h \|_E$ \\ \hline
		$8^3$   & $0.2500$  & $0.6233020$ & $0.148627$ & $4^3$ & $0.6434550$ & $0.044011$\\ 
		$16^3$ &  $0.1250$ & $0.6397600$ & $0.075046$ & $8^3$   & $0.6452138$ & $0.013348$  \\ 
		$32^3$ & $0.0625$ &  $0.6439755$ & $0.037636$ & $16^3$ & $0.6453773$ & $0.003826$  \\ 
		$64^3$ & $0.0313$ & $0.6450372$ & $0.018833$ & $32^3$  & $0.6453909$ & $0.001029$ \\
	       $128^3$ & $0.0156$ & $0.6453033$ & $0.009416$ & $64^3$ & $0.6453919$ & --- 
		\end{tabular}
	\end{center}
\end{table}

 \noindent We observe that the $\mathbb{Q}_1$  energy errors  are reducing by a factor of 2 
 with each grid refinement. This is consistent with the optimal $O(h)$ rate of 
 convergence predicted theoretically   in the case of   a $H^2$-regular problem.
 The $\mathbb{Q}_2$  energy errors  are reducing more rapidly with grid refinement.
 The observed rate is slightly less than  $O(h^2)$ which is the
 expected rate when solving  a $H^3$-regular   problem using triquadratic 
 approximation. 
 \end{examp}

\begin{examp}[Nonconvex domain]\label{ex.stairdomain}
A  $\mathbb{Q}_1$ finite element solution to  the Poisson problem (\ref{poisson})--(\ref{bc})  defined  on 
the domain  $D = {[-1,1]^3\setminus[-1,0)\times[-1,0) \times[-1,1] }$  
is shown in Fig.~\ref{fig.testproblem2}. In this computation
 the stair-shaped domain has been subdivided uniformly into  $32^3- 16^2\times32$ elements. The   dimension
 of the resulting linear algebra  system  is  27,489.  
 The MATLAB R2021b sparse direct solver solves this  
 system in about one fifth of a  second. The error plot illustrates the edge
 singularity in the solution along the reentrant corner edge $(x=0,y=0, z)$.

 \begin{figure}[!ht]
       \begin{center}
	\includegraphics[width=0.85\linewidth]{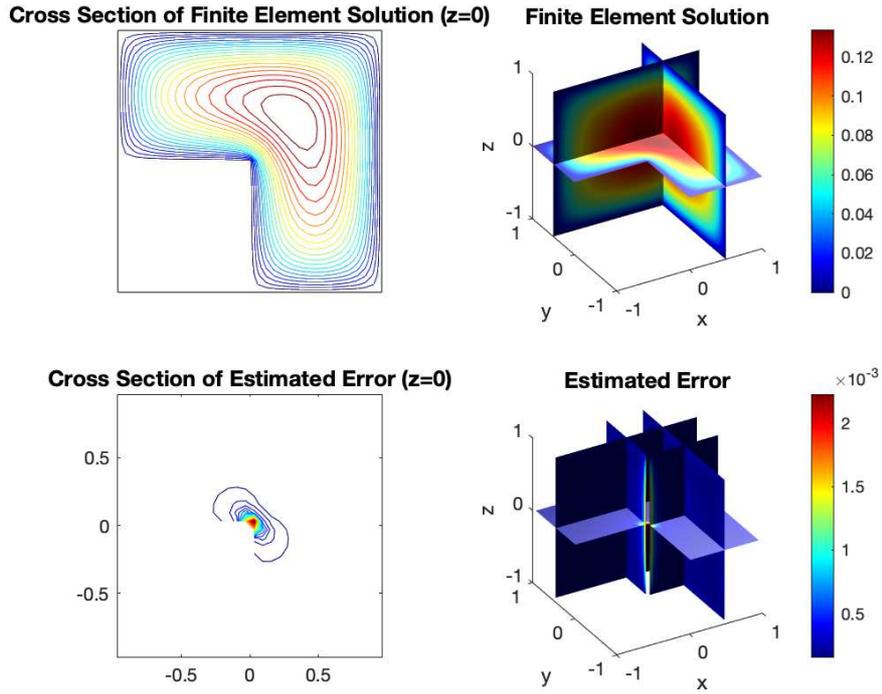}		 
	\end{center}\vspace{-4mm}
\caption{Solution and estimated energy error distribution for test problem~2.}\label{fig.testproblem2}
\end{figure}

  \begin{table}[!ht]
 \caption{Estimated energy errors for the staircase domain test problem.}
 \label{tab.testproblem2errors}
 \vspace{-2mm}
	\begin{center}
		\begin{tabular}{  c c c c | c c c} 
		& & \multicolumn{2}{c}{$\mathbb{Q}_1\qquad$} 
		&  \multicolumn{3}{c}{$\qquad \mathbb{Q}_2 $}   \\
		$n_e$ & $h$ & $\|u_h\|^2_E$ & $\| u_{ref}-u_h \|_E$  
		       & $n_e$ & $\|u_h\|^2_E$ & $\| u_{ref}-u_h \|_E$    \\ \hline
		$384$          & $0.2500$  & $0.2743216$ & $0.149663$ & $48$       & $0.2933030$ & $0.058461$ \\ 
		$3,072$      &  $0.1250$ & $0.2905480$ & $0.078566$  & $384$     & $0.2958987$ & $0.028670$  \\ 
		$24,576$    & $0.0625$ & $0.2949834$  & $0.041680$  & $3,072$  & $0.2964596$ & $0.016157$ \\ 
		$196,608$   & $0.0313$ & $0.2962188$ & $0.022402$ & $24,576$  & $0.2966318$ & $0.009424$ \\
		$1,572,864$& $0.0156$ & $0.2965759$ & $0.012030$ & $196,608$ &$0.2967206$& --- 
		\end{tabular}
	\end{center}
\end{table}

Estimated errors  for the second test problem are presented in Table~\ref{tab.testproblem2errors}.
 In this case we observe that the $\mathbb{Q}_1$  and $\mathbb{Q}_2$ energy errors  are both
 reducing by a factor of less than 2  with each grid refinement. This is  exactly what one
 would expect---the edge singularity limits the rate of convergence that is possible using uniform grids.  
 The second notable  feature is that the $\mathbb{Q}_2$ energy error is a factor of 2 smaller than 
 the $\mathbb{Q}_1$  error for the same number of degrees of freedom (the 
 results on the same horizontal line). This behaviour is also consistent with expectations; see
 Schwab~\cite{schwab98}.
  \end{examp}

 \begin{examp}[borehole  domain]\label{ex.boreholedomain}

A  $\mathbb{Q}_1$ finite element solution to  the Poisson problem (\ref{poisson})--(\ref{bc})  defined  on 
the cut domain  $D={[-1,1]^3\setminus (-\epsilon,\epsilon)\times[0,1]\times (-\epsilon,\epsilon)}$  
with $\epsilon=0.01$ 
is shown in Fig.~\ref{fig.testproblem3}. 

\begin{figure}[!ht]
       \begin{center}
	\includegraphics[width=0.85\linewidth]{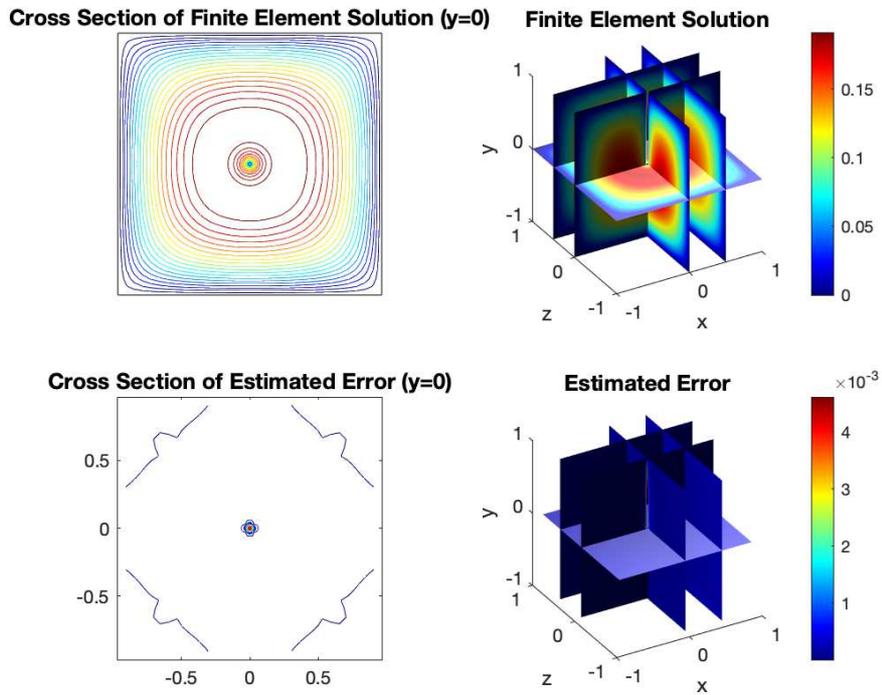}		 
	\end{center}\vspace{-4mm}
\caption{Solution and estimated energy error distribution for  test problem~3.}\label{fig.testproblem3}
\end{figure}

In this computation
 the borehole domain has been subdivided into  a tensor-product grid of  $48\times 32 \times 48$ elements
 with geometric stretching in the $x$ and $z$ direction so as  to  capture the geometry without using 
 an excessive number of elements. The grid spacing increases 
 from $h_{\min}=0.01=\epsilon$ within the hole to $h_{\max}=0.0625$ next to the boundary, so the maximum  
 element aspect ratio (adjacent to the hole) is 6.25. The  dimension  of the resulting linear algebra  system 
 is  85,833 and the MATLAB R2021b sparse direct solver  
solves the  system in about 6 seconds.   As anticipated, the error in the approximation is concentrated in 
a small region in the neighbourhood of the borehole, making  this a very challenging problem to solve efficiently.

\end{examp}
 
 \section{Structure of the software package}\label{sec.structure}
IFISS is designed for the MATLAB coding environment. This means 
that the source code is readable, portable and easy to modify. All local calculations (quadrature
in generating element matrices, application of essential boundary conditions, a posteriori error estimation) 
are vectorised over  elements---thus the code runs efficiently on contemporary Intel processor
architectures. IFISS3D has been developed for MATLAB (post 2016b) and tested with the 
current release (7.2) of Gnu Octave.  The main directory is called {\tt diffusion3D} and this needs to be added as a subdirectory of the main IFISS directory. The subdirectories of {\tt diffusion3D} are organised as follows.

$\bullet${\tt /grids/} \\
This directory contains all the functions associated with domain discretisation.  Three types of domain are included in the first release.  Introducing a new domain type is straightforward. A new function needs to be included  that saves nodal information (arrays  {\tt xyz}, {\tt bound3D}) and
(triquadratic) element information  ({\tt mv}, {\tt mbound3D}) in an
appropriately named datafile. This file will be subsequently read by an 
appropriate driver function associated with the specific PDE being solved.

\begin{figure}[!ht]
       \begin{center}
	\includegraphics[width=0.48\linewidth]{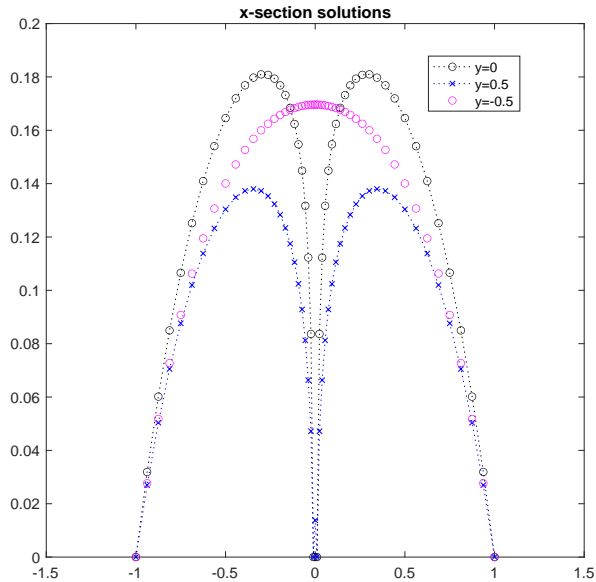}		 
	\end{center}
\caption{Incisions through the borehole test problem solution visualised in  Fig.~\ref{fig.testproblem3}.  The incisions  show the
solution in the  $x$ direction in the plane  $z=0$ for three different height values $y$.}\label{fig.incision3}
\vspace{-5mm} 
\end{figure}

$\bullet${\tt /graphs/} \\
This directory contains the functions associated with the visualisation  of the computed solution (nodal data) and 
the estimated errors (element data). The  tensor-product subdivision structure simplifies the code
structure substantially---plotting can be efficiently done using the built-in {\tt slice} functionality. Similarly, solution data
defined on a one-dimensional incision into the domain of interest can be plotted using the 
function {\tt xyzsectionplot}. An illustration is shown in Fig.~\ref{fig.incision3}.

$\bullet${\tt  /approximation/} \\
This directory contains all the functions associated with setting up the discrete matrix system associated with the 
PDE of interest. The functions {\tt femq1\_diff3D} and {\tt femq2\_diff3D} set up the stiffness and mass matrices
associated with the problems  discussed in Section~\ref{sec.referenceproblems}. Essential boundary
conditions are imposed by a subsequent call to the function {\tt nonzerobc3D}.   Extending
 the functionality by combining components of IFISS with IFISS3D to cover (a) nonisotropic diffusion and  (b) Stokes flow problems (using $\mathbb{Q}_2$--$\mathbb{Q}_1$
 mixed approximation) is a straightforward exercise. Efficient approximation of the solution
 of the heat equation on a three-dimensional domain can also be done with ease: 
 either using the  adaptive  time stepping functionality built into IFISS or using one
 of the ODE integrators  built into MATLAB. The functions associated with a posteriori error estimation 
can be found in four separate subdirectories associated with the four 
options described in Section~\ref{sec.specifics}.

$\bullet${\tt  /solvers/} \\
The  MATLAB sparse direct solver ($\backslash$)  has far from optimal complexity in a three-dimensional setting.  
This is explored in a case study in  the  next section.  
Algebraic multigrid (AMG) functionality is included in this directory to enable  exploration of an optimal solution strategy. 
If one does not have access to an efficient AMG setup routine then
the linear solver that is recommended  when the dimension of the system exceeds $10^5$ is 
MINRES (Minimum Residual) iteration preconditioned by an incomplete Cholesky factorisation of the system matrix with 
zero fill-in.\footnote{The 
incomplete factorisation function  {\tt ichol} provided in  MATLAB  R2021b is highly optimised.}.
 This strategy is encoded in the  {\tt it\_solve3D} function with a  residual stopping tolerance of $10^{-10}$.
 Solving the system in Example~\ref{ex.boreholedomain} using this strategy gives a solution in 66 iterations. The 
 associated CPU time is  half a second. This is over 10 times faster on a 2018 MacBook than the corresponding backslash solve! 
 
$\bullet${\tt  /test\_problems/} \\
This directory contains all the high-level driver functions such as {\tt diff3D\_testproblem} (the main driver). It also
contains the functions associated with the problem data (right-hand side function and essential
boundary specifications). The structure makes it straightforward to solve 
(\ref{poisson}) together with nonzero boundary data $u= g_D$ on $\partial D$. 

Help for IFISS is integrated into the MATLAB help facility. \rbl{Typing} {\tt helpme\_diff3D}
gives information on solving a Poisson problem in three dimensions. Starting from the main IFISS directory,
typing {\tt help diffusion3D/⟨subdirectory name⟩} gives a complete list of the files in that subdirectory.
\rbl{Using MATLAB, the function names are ``clickable'' to give additional information.}


\section{Case studies}\label{sec.casestudies}
Two important aspects of three-dimensional finite element approximation that can be investigated easily in IFISS3D are discussed in this section. 

\subsection{Effectivity of a posteriori error estimation strategies}

The effectiveness of hierarchical error estimation is well established in a two-dimensional setting; see, for
example,  Elman et al.~\cite[Table 1.4]{elman15}.  The IFISS3D software offers a choice of 4 such error estimation strategies in conjunction with  $\mathbb{Q}_1$ approximation. 
 Computed  error estimates  obtained when solving  the first test problem discussed in Section~\ref{sec.referenceproblems}
are presented in Table~\ref{tab.testproblem1err}. The four estimates are associated with the white
nodes shown in Fig.~\ref{fig:spacesYk}. The  estimated energy errors should be compared 
with the  reference energy errors  listed in Table~\ref{tab.testproblem1errors}.

\begin{table}[!ht]
 \caption{Computed error estimates $\eta_{\bullet}$ for the first test problem using four different hierarchical strategies.}
 \label{tab.testproblem1err}
 \vspace{-2mm}
 	\begin{center}
		\begin{tabular}{ c c c c c c } 
				$n_e$ & $h$ &$\eta_{\mathbb{Q}_2(h)} $&    $\eta_{\mathbb{Q}^r_2(h)} $ 
				&  $\eta_{\mathbb{Q}_1(h/2)}$
  & $\eta_{\mathbb{Q}^r_1(h/2)} $ \\
		\hline
			$8^3$& $0.1250$     & $0.150207$ & $0.137906$ &$0.129842$ &$0.115359$\\ 
			$16^3$& $0.0625$   & $0.075177$ & $0.069772$ & $0.065255$ &$0.058216$\\ 
			$32^3$& $0.0313$   &  $0.037648$ & $0.035050$ &$0.032655$ &$0.029215$\\ 
			$64^3$& $0.0156$   & $0.018837$ & $0.017561$  &$0.016326$ &$0.014631$\\
			$128^3$& $0.0078$ & $0.009420$ &$0.008789$   &$0.008161$ & $0.007321$
		\end{tabular}
	\end{center}
\end{table}

Table~\ref{tab.testproblem1eff} lists the associated effectivity indices.
 The indices get closer to 1 as the mesh is refined when the  $\mathbb{Q}_2(h)$ and 
 $\mathbb{Q}^r_2(h)$ strategies  are employed.
 On the other hand, the effectivity indices for the  $\mathbb{Q}_1(h/2)$ 
 and ${\mathbb{Q}^r_1(h/2)}$ strategies stagnate  around $0.87$ and $0.77$ respectively. All four error  estimates are correctly reducing by a factor of 2 with each grid refinement.
In light of these results, the relatively cheap ${\mathbb{Q}^r_2(h)}$ is set to be  the 
default option in IFISS3D.  Extensive testing on other problems
indicates that  this estimator consistently underestimates the error by a small amount.
 
\begin{table}[!ht]
 \caption{Effectivity indices $\theta_\bullet \coloneqq \eta_\bullet/  \| u_{ref}-u_h \|_E$ for test problem~1.}
 \label{tab.testproblem1eff}
 \vspace{-2mm}
	\begin{center}
		\begin{tabular}{  c c c c c } 
		$n_e$ & $\theta_{\mathbb{Q}_2(h)} $&  $\theta_{\mathbb{Q}^r_2(h)} $
			 &  $\theta_{\mathbb{Q}_1(h/2)}$& $\theta_{\mathbb{Q}^r_1(h/2)} $ \\
		\hline
		$8^3$ &$\boldsymbol{1.0106}$&$\boldsymbol{0.9279}$ & $\boldsymbol{0.8736}$&$\boldsymbol{0.7762}$\\ 
		$16^3$&$\boldsymbol{1.0017}$&$\boldsymbol{0.9297}$ & $\boldsymbol{0.8695}$&$\boldsymbol{0.7757}$\\ 
		$32^3$&$\boldsymbol{1.0003}$&$\boldsymbol{0.9313}$ & $\boldsymbol{0.8677}$&$\boldsymbol{0.7762}$\\ 
		$64^3$&$\boldsymbol{1.0002}$&$\boldsymbol{0.9325}$ & $\boldsymbol{0.8669}$&$\boldsymbol{0.7769}$		
		\end{tabular}
	\end{center}
\end{table}

All the errors reported in Table~\ref{tab.testproblem1err} were computed after making a
{\it boundary} element {\it correction}. This is a postprocessing step wherein the 
local problems associated with elements that have one or more boundary faces are 
modified so that  the (zero) error on the boundary is enforced as an essential boundary condition.
The motivation for making this correction is to recover the property of asymptotic
exactness in special cases.\footnote{An estimator is said to be {\it asymptotically exact} if 
the effectivity of the estimator tends to 1 when $h\to 0$.} The correction is, however, difficult to vectorise 
efficiently, raising the question as to whether it is worth including in a three-dimensional setting.

Computed effectivity indices for the special case of solving the Poisson problem 
\begin{align}\label{problem12}
		-\nabla^2 u &= f  \quad\text{in } D =[-1,1]^3
			\\
			u &= 0  \quad \text{on } \partial D \label{bc12}
	\end{align}
with the right-hand side function $f$  chosen so that the exact solution is the
triquadratic function
\begin{align}\label{solution12}
u(\boldsymbol{x}) &= (1- x^2)(1-y^2)(1-z^2),
\end{align}
are presented in Table~\ref{tab.testproblem12eff}.  The second and fourth columns are the 
results computed after making  the boundary correction. The asymptotic 
exactness of the $\mathbb{Q}_2(h)$ strategy can be clearly seen.
The third and fifth columns list the
results that are computed when the boundary correction is not made.  Comparing results with the 
second and fourth columns it is evident  that the boundary correction reduces the estimated error and, 
more importantly, that the size of  correction tends to zero in the limit $h\to 0$.  
The ${\mathbb{Q}^r_2(h)}$ strategy is not asymptotically exact so to speed up the 
computation  the default setting in IFISS3D is to simply {\it neglect} the boundary correction. Thus,
in the case of the finest grid in Table~\ref{tab.testproblem12eff}  (over 2 million elements) 
the  $\eta^*_{\mathbb{Q}^r_2(h)}$ error estimate is computed in
less than 9 seconds.   This is significantly less than the time taken to compute the
finite element solution  itself (the {\tt it\_solve3D} linear system solve  took over 23 seconds).

\begin{table}[!ht]
 \caption{Effectivity indices $\theta_\bullet \coloneqq \eta_\bullet/  \| u_{ref}-u_h \|_E$ for the Poisson problem
(\ref{problem12})--(\ref{solution12}).  The superscript $\theta^*_\bullet$ indicates that the boundary correction is
omitted in the computation of  $\eta_\bullet$.}
 \label{tab.testproblem12eff}
 \vspace{-2mm}
	\begin{center}
		\begin{tabular}{  c c c c c } 
		$n_e$ & $\theta_{\mathbb{Q}_2(h)} $& $\theta^*_{\mathbb{Q}_2(h)} $
			 &  $\theta_{\mathbb{Q}^r_2(h)} $& $\theta^*_{\mathbb{Q}^r_2(h)} $ \\
		\hline
		$16^3$  &$\boldsymbol{0.99944}$&$\boldsymbol{1.2914}$ & $\boldsymbol{0.97044}$& $\boldsymbol{0.99647}$\\ 
		$32^3$  &$\boldsymbol{0.99990}$&$\boldsymbol{1.1508}$ & $\boldsymbol{0.97087}$&$\boldsymbol{0.98289}$\\ 
		$64^3$  &$\boldsymbol{0.99998}$&$\boldsymbol{1.0771}$ & $\boldsymbol{0.97110}$&$\boldsymbol{0.97686}$\\
		$128^3$ &$\boldsymbol{1.00000}$&$\boldsymbol{1.0390}$ & $\boldsymbol{0.97123}$ &$\boldsymbol{0.97405}$	
		\end{tabular}
	\end{center}
\end{table}

\subsection{Fast linear algebra}

The solution of the (Galerkin) linear system is the computational bottleneck when solving a Poisson problem
in three dimensions. To illustrate this point, a representative timing comparison of the distinct solution 
components when solving the first test problem using  $\mathbb{Q}_1$ approximation with default settings
 is presented in Table~\ref{tab.timings}.  
The system assembly includes the grid generation. The overall time is the elapsed
time from start to finish and includes the time taken to visualise the solution and the associated error. 
What is immediately apparent is the fact that 
the system assembly times and the error estimation times scale approximately linearly with the
number of elements (or equivalently,  the  dimension  $n$ of the system matrix). The
backslash solve times, in contrast, grow like the square  of the number of elements on 
the most  refined grids. The memory requirement for the sparse factors of the system matrix  
also increases at a much faster rate than $O(n)$.  

\begin{table}[!ht]
 \caption{Representative component timings (in seconds) when solving test problem~1}
 \label{tab.timings}
 \vspace{-2mm}
	\begin{center}
		\begin{tabular}{  r c c c c } 
		$n_e$ & assembly &  solve  ($\backslash$)
			 &  estimation & overall  \\
		\hline
		$16^3$    &$0.171$ & $\boldsymbol{0.016}$ & $0.617$ & $\boldsymbol{1.10}$\\ 
		$32^3$    &$0.669$ & $\boldsymbol{0.339}$ & $1.481$ & $\boldsymbol{4.81}$\\ 
		$64^3$    &$7.100$&$\boldsymbol{25.29}$ & $11.49$&$\boldsymbol{58.1}$\\
		$128^3$  &$90.16$  &$>\boldsymbol{10^3}$ & $98.45$ & $>\boldsymbol{10^3}$	
		\end{tabular}
		
	\end{center}
\end{table}

The optimal  $O(n)$ complexity of the overall solution algorithm can be recovered by solving 
the linear system using a short-term Krylov subspace iteration such as MINRES in combination with an
algebraic multigrid (AMG) preconditioning strategy; see~\cite[sec.\,2.5.3]{elman15}. 
The  set up phase of  AMG is a recursive procedure: heuristics associated with  algebraic relations 
(``strength of connections'') between the unknowns are used to generate a sequence of
progressively coarser representations  $A_{\ell}$, $\ell =1,2,\ldots, L$ of the 
Galerkin system matrix $A \coloneq A_1$.  The solution (preconditioning)  phase approximates 
the action on the inverse of $A$ on a vector by cycling through the associated grid sequence.
At each level, a fixed-point iteration (typically point Gauss-Seidel) is applied
to ``smooth'' the residual  error that is generated by interpolation or restriction of the error
vector generated at the previous level. If coarsening is sufficiently rapid then the 
work associated with the preconditioning step will be proportional to the \rbl{number
of unknowns.} 

 The algorithmic complexity of any AMG coarsening strategy can be characterised by 
 a few parameters. First, the  {\it grid complexity} is defined as 
$$c_G \coloneq \frac{1}{n_{1}} \sum_{\ell=1}^L n_{\ell},$$
where   $n_{\ell}$ is the dimension of the coarse grid matrix $A_{\ell}$ at level $\ell$.
Starting from a uniform grid of  $\mathbb{Q}_{1}$ elements, if full coarsening 
(in each spatial direction) is done at each level then $n_{\ell}$ reduces by a factor of $8$ at each level, in which 
case we obtain $c_{G} \approx {8}/{7}$.  A value of $c_{G}$
higher than this suggests that coarsening has not been done isotropically.
The {\it operator complexity} is typically defined by
$$c_A \coloneq \frac{1}{{\rm nnz}(A_{1})} \sum_{\ell =1}^{L} {\rm nnz}(A_{\ell}),$$
where ${\rm nnz}(A)$ is the number of nonzeros in the matrix.
This parameter provides information about the associated storage requirements for the coarse grid 
matrices generated. If uniform coarsening is done and the coarse grid 
matrices correspond to the usual finite element discretisation on those grids then we would 
expect $c_{A} \approx c_{G}$. In practice however, the coarse grid matrices become progressively
denser, with larger stencil sizes, as the level number increases. If the matrices become 
too dense then this may cause an issue with the computational cost of  applying the
 chosen smoother.  To  quantify this, the  {\it average stencil size} 
$$c_S: = \frac{1}{L} \sum_{\ell=1}^L \frac{{\rm nnz} (A_{\ell})}{n_{\ell}},$$
should be compared with the average stencil size at the finest level, that is,
$$c_1  \coloneq \frac{{\rm nnz}(A_{1})}{n_{1}}. $$

An implementation of  the coarsening strategy developed by 
Ruge and St{\"u}ben~\cite{ruge87}  is included in the 
\rbl{IFISS software}.
The corresponding IFISS3D function {\tt amg\_grids\_setup3D} 
can be edited in order to explore algorithmic options 
or change the default  threshold parameters.  
The MATLAB implementation of the coarsening algorithm  is far from optimal however. 
There is a marked deterioration in performance when solving problems  on fine grids 
that is evident even when solving Poisson problems in two dimensions.  
 To address this issue an interface to the efficient 
Fortran~95  implementation~\cite{hsl} of the same coarsening algorithm
is included in IFISS3D.\footnote{The
HSL\_MI20 source code and associated MATLAB interface is freely available to 
staff and students of recognised educational institutions. The inclusion of the
compiled {\tt mex} file is prohibited by the terms of the HSL academic licence.}

\vspace{-2mm}
\begin{table}[ht!]
\caption{AMG grid complexity data and representative linear solver timings  for
test problem~1.}
\vspace{-2mm}
\label{tab.amgresults1}
	\begin{center}
\begin{tabular}{ c r r r r | r r r } 
& \multicolumn{4}{c}{$\mathbb{Q}_1$} 
		&  \multicolumn{3}{c}{$\mathbb{Q}_2 $}  \\
  $n_e$ & $16^3$  & $32^3 $  & $64^3$  &  $128^3 $  &  $16^3 $  & $32^3$  &  $64^3 $ \\ \hline 
 $n_{1}$ &  4,913  & 35,937 & 274,625 & 2,146,689  & 35,937 & 274,625 & 2,146,689\\\
 $c_{1}$  & 16.49 & 21.14 &  23.90 & 25.41 & 47.06 &  54.96 & 59.33 \\ \hline 
  $L$  &  7  &  11  &  16 & 14 &  11 & 14 & 16 \\
  $c_{G}$  &  1.24  &  1.32 &  1.37 & 1.39  & 1.27 & 1.32 & 1.32\\
  $c_{A}$  &  1.58 &  1.77  & 1.88  & 1.93  &1.61 & 1.70 &1.76 \\
  $c_{S}$ &  22.9 &  39.2 &  62.8   & 97.6 & 65.59 & 106.1 & 138.8 \\
   setup${}^*$  time   & 0.01 & 0.07 &  0.58 & 5.12 & 0.15 & 1.22 & 11.49\\ \hline
    total  time   & \textbf{0.02}  & \textbf{0.16} &    \textbf{1.49} &  \textbf{13.07} 
                       & \textbf{0.42} &    \textbf{3.48} &  \textbf{33.60}
\end{tabular}
\end{center}
\end{table}

AMG complexities and timings obtained when solving test problem~1 
are presented in Table~\ref{tab.amgresults1}. The total times reported for the $\mathbb{Q}_1$ approximation
 should be compared  with the corresponding  backslash solve timings recorded in Table~\ref{tab.timings}.
 The setup times were recorded  using the interface to the  HSL\_MI20 code and
 scale  close to linearly with the problem dimension. This behaviour is consistent with the results
that  were reported for the same problem and discretisation
by Boyle et al.~\cite[Ex.\,4.5.1]{hslcode}.  Looking at the AMG grid data in Table~\ref{tab.amgresults1},
we note that the grid complexity is under control and stays close to  the optimal value using
either of the two approximation strategies. 
The operator complexity results are also  encouraging---increasing slowly as the dimension 
of the problem is increased. The growth in the average stencil size is not unexpected.
 The value of $c_S$ is within a factor of  3 of $c_1$ when using $\mathbb{Q}_2 $ approximation
 on the finest grid.  The difference  between the highlighted total time and the associated setup time
 is  the time taken by preconditioned MINRES to reach the 
 residual stopping  tolerance of $10^{-10}$ with the  default smoothing parameters 
 (that is, two pre- and post-smoothing sweeps using point Gauss-Seidel).
 
 The (Ruge--St{\"u}ben) AMG  coarsening strategy is designed to be effective even in cases where the discretised problem has strong anisotropy in the grid spacing.
 This is validated  by  the encouraging  timing results for the borehole problem that are presented in Table~\ref{tab.amgresults3}. Here, $l$ denotes the level of refinement of the finite element mesh.

\vspace{-2mm}
 \begin{table}[ht!]
\caption{AMG grid complexity data and representative linear solver timings  for  
test problem~3.}
\vspace{-2mm}
\label{tab.amgresults3}
	\begin{center}
\begin{tabular}{ c r r r r} 
& \multicolumn{4}{c}{$\mathbb{Q}_1$} \\
               & $l=2$  & $l=3 $  & $l=4$  &  $l=5$\\ \hline 
 $n_{1}$ &   85,833   & 365,625  &  1,264,329  &  3,888,153 \\\
 $c_{1}$  & 22.41 & 24.15 & 25.08  & 25.61 \\ \hline 
  $L$  &  15  &  16  &  18 & 18 \\
   setup${}^*$  time   & 0.44 & 2.42 &  8.48 & 26.47 \\ \hline
    total  time   & \textbf{1.08}  & \textbf{6.89} &    \textbf{24.05} &  \textbf{79.20}
\end{tabular}
\end{center}
\end{table}

\newpage
\section{Summary and future developments}\label{sec.summary}

The IFISS3D toolbox extends the capabilities of IFISS \cite{ifiss} to facilitate the numerical solution of elliptic PDEs on three-dimensional spatial domains that can be partitioned into hexahedra. In particular, it allows users to investigate the convergence properties of trilinear ($\mathbb{Q}_{1}$) and triquadratic ($\mathbb{Q}_{2}$) finite element approximation for test problems whose solutions have varying levels of spatial regularity and the performance of a range of iterative solution algorithms for the associated discrete systems, including an optimal AMG solver. For $\mathbb{Q}_{1}$ approximation, the effectivity of four distinct state-of-the-art hierarchical error estimation schemes can also be explored. The IFISS3D software is structured in such a way that, when integrated into the existing IFISS software, users can easily solve a range of other PDE problems, including time-dependent ones, using $\mathbb{Q}_{1}$ and $\mathbb{Q}_{2}$ elements on three-dimensional spatial domains. IFISS together with IFISS3D is intended to be useful as a teaching tool, and can be used to produce matrices of arbitrarily large dimension for testing linear algebra algorithms. Future developments of IFISS3D will be documented on \rbl{GitHub}.

\bibliographystyle{siam}
\bibliography{references}

\end{document}